\documentclass[numreferences]{kluwer}

\newtheorem{theorem}{Theorem}
\newtheorem{lemma}[theorem]{Lemma}
\newtheorem{remark}[theorem]{Remark}

\newtheorem{proposition}[theorem]{Proposition}
\newtheorem{corollary}[theorem]{Corollary}

\newcommand\pr{\mathbb{P}^n}
\newcommand\prm{\mathbb{P}^m}
\newcommand\pri{\mathbb{P}^1}
\newcommand\cO{\mathcal{O}}

\newcommand\Op{O_p}
\newcommand\Cp{\overline{\mathbb{Q}}_p}
\newcommand\Fp{\mathbb{F}_p}
\newcommand\oFp{\overline{\mathbb{F}}_p}

\newcommand\Qp{\mathbb{Q}_p}

\newcommand\OpT{\Op[T_0,\ldots,T_n]}

\newcommand\cp{\mathfrak{p}}

\newcommand\cq{\mathfrak{q}}
\newcommand\cm{\mathfrak{m}}
\newcommand\oC{\overline{C}}

\begin{document}
\begin{article}
\begin{opening}


\title{Periodic points for good reduction maps on
curves\thanks{The first author was supported
by EPSRC Grant GR/M 49588, LMS Scheme 4 grant 4627
and FWF Project P14379; the authors
thank two anonymous referees for their comments.}}
\author{Manfred \surname{Einsiedler}
\email{manfred.einsiedler@univie.ac.at}}
\institute{Institut f\"ur Mathematik,
  Universit\"at Wien, Strudlhofgasse 4, A-1090 Vienna, Austria}
\author{Graham \surname{Everest}\email{g.everest@uea.ac.uk}}
\author{Thomas \surname{Ward}\email{t.ward@uea.ac.uk}}
\institute{School of Mathematics, University of East
  Anglia, Norwich NR4 7TJ, United Kingdom}


\classification{Mathematics Subject Classification (2000)}{11S05, 37F10}

\begin{abstract}
The periodic points of a morphism of good reduction for a smooth
projective curve with good reduction over $\Cp$ form a discrete
set. This is used to give an interpretation of the morphic height
in terms of asymptotic properties of periodic points, and a
morphic analogue of Jensen's formula.
\end{abstract}
\keywords{Morphism of good reduction, Periodic points, Morphic
heights}
\end{opening}

\section{Introduction and Results}\label{intro}

In this paper we study the behaviour of periodic points for a
morphism on a smooth curve over $\Cp$,
the algebraic closure of $\Qp$, with
$p$-adic norm $|\cdot|$ normalized to have
$\vert p\vert=1/p$.
Under a
regularity condition, we prove that the asymptotic distance of a
given point to the periodic points is equal to one in
a suitable metric.
This result generalises the case of polynomial
morphisms on the projective line in~\cite{nynts}.

Let $O_p=\{z\in\Cp\mid\vert
z\vert\le1\}$ be the ring of integers in $\Cp$, with maximal ideal
$\cp=\{z\in\Cp\mid\vert z\vert<1\}$. Identify the quotient
$O_p/\cp$ with the algebraic closure $\oFp$ of $\Fp={\mathbb
Z}/p{\mathbb Z}$. Let $\pr(\Cp)$ denote $n$-dimensional projective
space over $\Cp$, and write $X=(X_0,\dots,X_n)\in\Cp^{n+1}$ for
the homogeneous coordinates of a point $[X]\in\pr(\Cp)$.
It will be useful always to choose
the homogeneous coordinates to have
\begin{equation}\label{normalizechoiceofX}
\vert X\vert=\max\{\vert X_0\vert,\dots,\vert
X_n\vert\}=1.
\end{equation}
Writing $x=[X],y=[Y]$, define a function on $\pr(\Cp)$ by
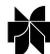
\begin{equation}\label{definemetric}
\Delta(x,y)=\max_{i,j}|X_iY_j-X_jY_i|.
\end{equation}
Then, just as in the one-dimensional
case (see \cite{call-silverman-1993} for example), $\Delta$ is a
metric. In this metric,
projective space $\pr(\Cp)$ has diameter one.
The regularity property used here -- {\sl good reduction} --
will be defined in the next section.

A point $x\in C$ is a {\sl fixed point for} a map
$\phi:C\to C$ if $\phi(x)=x$, is a {\sl point of period
$n$} if $\phi^n(x)=x$, and has {\sl least period} $n$
if the orbit $\{x,\phi(x),\phi^2(x),\dots\}$ has
cardinality $n$.

\begin{theorem}\label{mainth}
 Let $C\subseteq\pr(\Cp)$ be an irreducible smooth curve with good reduction,
 and let $\phi:C\rightarrow C$ be a
 morphism of good reduction with degree $d>1$.
Then for every point $x\in C$ there is a constant
$\kappa_x$ with
$$
\Delta(x,y)\ge\left(\kappa_x\vert n\vert\right)^{1/n}
$$
for any $\phi$-periodic point $y$ with least period $n$.
\end{theorem}

\begin{corollary}
For every point $x\in C$, and every $r<1$, the number
of points $y\in C$ that are periodic under $\phi$ and
satisfy $\Delta(x,y)<r$ is finite.
\end{corollary}

The Riemann--Hurwitz formula~\cite[Sect.~IV.2]{hart}
shows that
$$
(2-2g)(d-1)\ge0
$$
where $g$ is the genus of the curve $C$,
since the
ramification divisor of $\phi$ is non-negative.
In particular,
the hypothesis of Theorem~\ref{mainth} implies that
$C$ has genus $0$ or $1$.

In the case of the projective line $C={\mathbb
P}^1(\Cp)$ there are infinitely many periodic points (see
\cite[Theorem 6.2.2]{beardon-1991}; it is enough to
prove this over any field). The proof of Theorem
\ref{mainth} simplifies in this case.
If the point
at infinity of the projective line is fixed by $\phi$,
then the periodic points of the morphism are related to the
{\sl morphic height}
associated to $\phi$ (cf. Theorem \ref{rational_case} below; the morphic
height is defined in~\cite{call-silverman-1993}).

A rational function $\phi$ has {\sl good
reduction} modulo $p$ if it can be written in homogeneous
coordinates in the form
\[
 \phi(X,Y)=(F(X,Y),G(X,Y))
\mbox{ with } F,G\in\Op[X,Y]
\]
where $F$ and $G$ have no common root modulo $\cp$ (see
\cite{benedetto}, \cite{call-silverman-1993} and
\cite{thesis}). Write
$\lambda_{\phi,p}$ for the local morphic height
(sometimes called the canonical local height in
the literature) as defined in
\cite{call-silverman-1993}.

\begin{theorem}\label{rational_case}
 Let $\phi:{\mathbb P}^1(\Cp)\rightarrow{\mathbb P}^1(\Cp)$ be a
 rational function of good reduction and degree $d>1$. Assume that the
 point $(0,1)$ at infinity is fixed under $\phi$. Fix a
 point $x\in\Cp$, and a sequence $y_n\in\Cp$ of points of least period
 $n$ under $\phi$. Then
 \[
  \frac{|x-y_n|}{\max(1,|x|)\max(1,|y_n|)}\rightarrow 1\mbox{ for
  }n\rightarrow\infty.
 \]
 Moreover,
 \[
  \log|x-y_n|\rightarrow\lambda_{\phi,p}(x)=\log^+|x|_p
 \]
 and, if $x$ is not a periodic point, then
 \[
 \frac{1}{d^n}\log|f_n(x)-xg_n(x)|\rightarrow\lambda_{\phi,p}(x)=\log^+|x|
 \]
 where $\phi^n=\frac{f_n}{g_n}$, so $f_n(t)-tg_n(t)$ is the
 polynomial whose roots are exactly the periodic points of period $n$.
\end{theorem}

Notice that there are infinitely many points
whose least period exceeds $1$ (cf. Remark~\ref{remperpts}).
Theorem~\ref{rational_case} is a morphic analogue of Jensen's
classical formula (see Section~\ref{examplessection}).
Finally we should point out that a paper
of Lubin~\cite{lubin} contains results closely
related to those presented here, and
Hua--Chieh Li has made an extensive study~\cite{li1},
\cite{li2}, \cite{li3},
of periodic points for $p$-adic power series,
mainly aimed at counting the points of given
period.

\section{Proofs of theorems}\label{proof_main}

In this section we prove Theorems \ref{mainth} and
\ref{rational_case} assuming some results on good reduction curves
and uniformizers that will be proved later. Recall that
$X=(X_1,\ldots,X_n)$ always denotes the homogeneous coordinates of a
point $x=[X]\in C$ chosen so that Equation (\ref{normalizechoiceofX})
holds.

Let $\pi:\pr(\Cp)\to\pr(\oFp)$ be the reduction map,
defined by
\[
\pi(x)=[X_0+\cp,X_1+\cp,\dots,X_n+\cp],
\]
which is well-defined by (\ref{normalizechoiceofX}). Let $C$ be an
irreducible projective curve in $\pr(\Cp)$, with ideal of
relations $I=I(C)$. Let $J=I\cap\OpT$ be generated by the forms
$f_1,\dots,f_t$. Fix a point $y\in C$, and assume without loss of
generality that $Y_0\neq0$. The curve is {\sl non-singular at $y$}
if \begin{equation}\label{rankcondition}
\mbox{rank}\left(\frac{\partial g_i}{\partial
U_j}(y)\right)_{i,j}=n-1,
\end{equation}
where $g_i(U_1,\dots,U_n)=f_i(1,U_1,\dots,U_n)$. The curve $C$ is
{\sl smooth} if it is non-singular at every point.

Define $\overline{J}=J\mbox{ mod }\cp\subset\oFp[T_0,\ldots,T_n]$,
and let $\overline{C}\subset\pr(\oFp)$ be the variety defined by
the ideal $\overline{J}$ (which of course may not in general
coincide with the ideal of
relations of the algebraic set $\overline{C}$). The curve $C$ has
{\sl good reduction} if (\ref{rankcondition}) holds mod $\cp$ for
every $\bar{y}\in\overline{C}$. From now on we will assume that
the curve $C$ is smooth with good reduction.

The metric (\ref{definemetric}) on an integral affine piece
simplifies
as follows: for $x=[X],y=[Y]\in\pr(\Cp)$ with $X_0=1, Y_0=1$ and
$X_i,Y_i\in O_p$,
\[
 \Delta(x,y)=\max_i|Y_i-X_i|.
\]

A rational function $f={F}/{G}$ on $C$ is defined by two forms
$F$ and $G$ in $\Cp[T_0,\ldots,T_n]$ of the same degree with
$G\notin I(C)$. A rational function is {\sl regular at} $x\in C$
if there are two forms $F',G'$ with $G'(X)\neq 0$ and $FG'-F'G\in
I$ (in other words $f={F'}/{G'})$. Moreover, $f$ is {\em
regular at} $\bar{x}\in\oC\subseteq\pr(\oFp)$ if there are two
forms $F',G'\in\Op[T_0,\ldots,T_n]$ such that $G'(\bar{X})\neq 0$,
$FG'-F'G\in I$, and $\bar{X}\in\oFp^{n+1}$ is a homogeneous
coordinate of $\bar{x}$. In that case
$\overline{f}={\overline{F}}/{\overline{G}}$ defines a
rational function on $\oC$ which is regular at $\bar{x}\in\oC$.

A {\sl uniformizer} of $C$ at $x\in C$ is a rational function $z$
which is regular at $x$, such that the vector $(\frac{\partial
z}{\partial U_i})_i$ is not in the image of the matrix in
(\ref{rankcondition}).

\bigskip\noindent{PROPOSITION \ref{bijective}} (cf. Section
\ref{Sect_unif}.) {\it
 Fix a point $x\in C$.
 Then there exists a uniformizer $z$ of $C$ at $x$ such that
 $z$ is regular at $\pi(x)$ and $\overline{z}$ is a uniformizer at
 $\pi(x)$.
 The restriction $z:U\to\cp$ of $z$ to $U=\{y\in C\mid
 \Delta(y,x)<1\}$
 is a bijection and its inverse, in each affine coordinate, is
 a convergent power
 series with coefficients in $O_p$. Hence $\Delta(y,y')=|z(y)-z(y')|$ for
 all $y,y'\in U$.}
\smallskip

For example, if the curve is $\pri(\Cp)$ and $X=(1,0)$, one may choose
$z(y)=\frac{Y_1}{Y_0}$ for the uniformizer in Proposition~\ref{bijective}.

Now let $\phi:C\rightarrow C$ be a morphism, defined by forms
\[
 (F_0,\ldots,F_n)\in\Cp[T_0,\ldots,T_n]^{n+1}
\]
of the same degree, with $F_j\notin I(C)$ for some $j$. Since
$\phi$ is a map from $C$ to $C$,
$f(F_0,\ldots,F_n)\in I(C)$ for every $f\in I(C)$. In order for
$\phi$ to be defined on all of $C$, at each $x=[X]\in C$ there must be
a representation
\begin{equation}\label{representation}
 (G_0,\ldots,G_n)\in\Cp[T_0,\ldots,T_n]^{n+1}
\end{equation}
of the morphism, with
\begin{equation}\label{represents_the_same}
 F_iG_j-F_jG_i\in I(C)\mbox{ for }0\leq i,j\leq n
\end{equation}
and
\begin{equation}\label{is_defined_at_x}
 (G_0(X),\ldots,G_n(X))\neq (0,\ldots,0).
\end{equation}
Here~(\ref{represents_the_same}) means that the forms
$(G_0,\ldots,G_n)$ define the same map as the forms
$(F_0,\ldots,F_n)$ do, and (\ref{is_defined_at_x}) means that
$\phi$ is well defined at $x$.

The morphism $\phi:C\rightarrow C$ has {\sl good reduction} if for
every $\bar{x}\in\oC$, there is a representation
\[
 (G_0,\ldots,G_n)\in\Op[T_0,\ldots,T_n]^{n+1}
\]
of $\phi$ satisfying (\ref{represents_the_same}) such that
\[
 (G_0(\bar{X}),\ldots,G_n(\bar{X}))\notin\cp^{n+1}.
\]
Here $\cp^{n+1}$ denotes the
$(n+1)$-fold Cartesian product of
$\cp$ and $\bar{X}$ is a homogeneous coordinate for
$\bar{x}\in\oC\subseteq\oFp^{n+1}$.

In the case of $C=\pri(\Cp)$, there exists a canonical
representation $(F_0,F_1)$ satisfying (\ref{is_defined_at_x}) for
all $x\in\pri(\Cp)$. We can assume that $F_0,F_1\in\Op[T_0,T_1]$
and at least one of the two polynomials has a coefficient in
$\Op^{\times}$. Then the morphism $\phi$ has good reduction if and only
if the two forms $\overline{F_0},\overline{F_1}$ do not have a
common zero on $\pri(\oFp)$ -- they define a rational function on
$\pri(\oFp)$.

\begin{remark}
 A morphism $\phi:C\rightarrow C$ has good reduction in the above sense if and only
 if it extends to a morphism over the scheme
 $\mbox{Spec}(O_p)$.
\end{remark}

From now on we will assume that $\phi$ is a morphism of good
reduction.

Let $K$ denote the field of rational functions on $C$. The degree
of the morphism $\phi$ is defined as the degree of the field
extension $[K:\phi^*(K)]$, where
\begin{eqnarray*}
 \phi^*&:&K\rightarrow K\\
  &&f\mapsto f\circ\phi
\end{eqnarray*}
is the map induced by $\phi$. Alternatively, one can define the
degree as the common number of pre-images (counted with
multiplicities) of points under the map $\phi$.

In the case of a rational function $\phi$ on $\pri(\Cp)$ this is
again obvious since the degree $d$ of $\phi$ coincides with the
degree of the forms in the canonical representation $(F_0,F_1)$.

\bigskip\noindent{COROLLARY \ref{uni_power}} (cf. Section
\ref{Sect_unif}.) {\it
 Let $\phi:C\rightarrow C$ be a morphism of good reduction, and let
 $x\in C$ be a fixed point. Then there exists a
 uniformizer $z$, also satisfying Proposition \ref{bijective},
 such that
 \[
  \phi^*z=z\circ\phi=\sum_{i=1}^ {\infty}a_iz^i\mbox{ with
  }a_i\in\Op.
 \]
}\smallskip

Again if $C=\pri(\Cp)$ and $\phi$ is a rational function on
$\pri(\Cp)$ with good reduction, it is easy to see that Corollary
\ref{uni_power} holds: Assume $X=(1,0)$, then
$z(y)=\frac{Y_1}{Y_0}$ is a uniformizer at $x$, and
\[
 \phi^*z(y)=\frac{F_1(Y)}{F_0(Y)}.
\]
Since $x$ is a fixed point, $F_1(X)=0$ and $F_0(X)\neq
0$. Moreover, $\phi$ has good reduction and therefore
$F_0(X)\notin\cp$. Changing to the affine coordinate $z$,
write $\phi^*z=\frac{f_1(z)}{f_0(z)}$ with $f_1(0)=0$ and
$|f_0(0)|=1$. Since the constant term of $f_0$ is a unit in $\Op$,
the polynomial $f_0$ is a unit in $\Op[[z]]$. So the rational
function $\phi^*z\in\Op[[z]]$ satisfies Corollary~\ref{uni_power}.

\begin{lemma}\label{localphi_neu}
 Let $f(z)=\sum_{i=1}^\infty a_iz^i$ be a power series with
 coefficients $a_i\in \Op$. Then, for any $n>1$,
 \[
  f^n(z)=a_1^nz+z^2g_n,
 \]
 where $g_n\in\Op[[z]]$.
 Assume now that $a_1=1$ and let $e$ be the first index $e>1$ with
 $a_e\neq 0$, so that $f(z)=z+z^eg$ for some $g\in\Op[[z]]$.
 For any $n>1$ there exists a power series $h_n\in\Op[[z]]$ with
 \[
  f^n(z)=z+nz^eg+z^{2e-1}h_n.
 \]
\end{lemma}

This may be seen by a simple induction argument.
For the proof of the second statement, notice
that $g(z+z^eF(z))=g(z)+z^eG(z)$ for
some $G$ depending on $g$ and $F$.

\begin{remark}\label{remperpts} Some
rather general properties of periodic points
in the setting are needed later and assembled
here.
\begin{enumerate}
\item The second statement in Lemma~\ref{localphi_neu}
can be used to show that if a point $x$ is
a multiple root of $f^n(z)-z=0$, then
it is never a higher multiplicity root of
any other equation of the form $f^m(z)-z=0$
(notice this is only true for multiplicity two
or higher). In particular, if $x\in C$ is a point
of period $n$ for $\phi$ with multiplicity two or higher,
then that multiplicity cannot increase
when $x$ is viewed as a point with period $m>n$.
\item If the curve is ${\mathbb
P}^1(\Cp)$ then many of the roots of
$\phi^n(z)-z=0$ are genuine points with least
period $n$, and in particular there are infinitely
many points whose least period exceeds any given number.
\item Any map for which the number of points of period
$n$ grows exponentially fast will have the
same exponential rate of growth in the number of
points of least period $n$ (cf.~\cite{pw}).
\end{enumerate}
\end{remark}

\begin{lemma}\label{referee}
 Let $f(z)=\sum_{i=1}^\infty a_iz^i$ be a power series with
 coefficients $a_i\in \Op$ such that $|a_1-1|<p^{-1/(p-1)}$. Then every
 periodic point $y\in\cp\setminus\{0\}$ of $f$ of least period $n\geq
 1$ satisfies
 \[
  |y|\geq (\kappa |n|)^{{1}/{n}},
 \]
 where the constant $\kappa>0$ does not depend on $n$.
\end{lemma}

\begin{pf}
 The assumptions on $f$ imply that $|f(y)|\leq |y|$
 for any $y\in\cp$. So if $y\in\cp$ is a periodic point of least
 period $n$, then the points
$y_1=y, y_2=f(y),\ldots, y_n=f^n(y)$ along the orbit of $y$ all have
 the same norm.

 Suppose first that $a_1\neq 1$. From the $p$-adic
logarithm (see~\cite[Sect.~IV.1]{koblitz}) it follows that
 for every integer $n\geq 1$, $|a_1^n-1|\geq\kappa |n|$,
where $\kappa=|\log_p a_1|$.
 Consider the power series
 \begin{equation}\label{dustman}
  F_n(z)=f^n(z)-z=\sum_{i=1}^\infty b_iz^i\in\Op[[z]];
 \end{equation}
the first nontrivial term for this series is $b_1=a_1^n-1$.
If $y$ is a periodic point in $\cp\setminus\{0\}$ of least period $n$,
then all the points
 $y_i$ on the orbit of $y$ are roots of the equation $F_n(z)=0$.
 From the usual Newton
polygon arguments (see for example~\cite[Sect.~IV.4]{koblitz})
we see that $\log |y|=\log|y_i|<0$ equals
 one of the slopes of the Newton polygon of $F$, say the slope between
 the points $P_k$ and $P_\ell$ defined by the coefficients $b_{k}$
and $b_\ell$ with $k<\ell$ (cf. Figure 1).
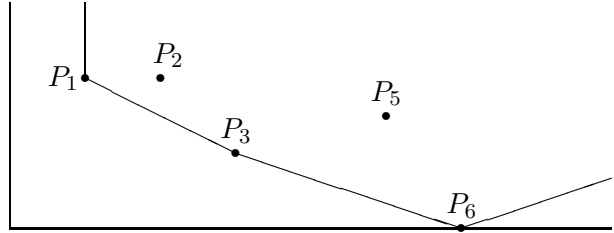
\begin{figure}[h]
\begin{center}
\setlength{\unitlength}{1mm}
\begin{picture}(80,30)
\put(0,0){\line(1,0){80}}
\put(0,0){\line(0,1){30}}
\put(10,30){\line(0,-1){10}}
\put(10,20){\circle*{1}}
\put(10,20){\line(2,-1){20}}
\put(30,10){\circle*{1}}
\put(30,10){\line(3,-1){30}}
\put(60,0){\circle*{1}}
\put(60,0){\line(3,1){20}}
\put(5,19){$P_1$}
\put(19,22){$P_2$}\put(20,20){\circle*{1}}
\put(28,12){$P_3$}\put(48,17){$P_5$}\put(50,15){\circle*{1}}
\put(58,2){$P_6$}
\end{picture}
\end{center}
\caption{\label{f:newton} The Newton polygon of $F$ is defined by
the points $P_i=(i,-\log|b_i|)$. The slopes determine the norms of
the roots.}
\end{figure}
Since $\ell-k$ is exactly the
number of roots whose norm equals the slope, we get
 $\ell-k\geq n$. Furthermore $|b_1|\leq |b_k|$, because
 otherwise the point defined by $b_k$ is higher
 in the Newton polygon and the slope between the points
 $P_k$ and $P_\ell$ is positive -- contradicting $|y|<1$.
 Furthermore, $|b_\ell|\leq 1$ and
 \[
  \frac{\log|b_1|}{n}\leq\frac{\log|b_1|-0}{\ell-k}
\leq\frac{\log|b_k|-\log|b_{\ell}|}{\ell-k}=\log|y|.
 \]
Together with the estimate on $b_1=a_1^n-1$, this proves the lemma
in the case $a_1\neq1$.

If $a_1=1$, choose
$\kappa=\vert b_e\vert$, where
$b_e$ is the first non-zero
coefficient in~(\ref{dustman}),
use the second equation
in Lemma~\ref{localphi_neu}, and deduce in a similar fashion the same
inequality.
\end{pf}

We are grateful to a referee for pointing out that
the condition
$$
|a_1-1|<p^{-1/(p-1)}
$$
in the hypotheses
of~Lemma~\ref{referee} can be removed entirely
by arguing as follows.
Let $m\ge1$ be the smallest integer with
$|a_1^m-1|<p^{-1/(p-1)}$. If $a_1^m\neq1$,
then choose $\kappa=\vert\log_p(a_1^m)\vert$.
For $n\ge1$,
$$
\vert a_1^n-1\vert\ge\vert a_1^{mn/d}-1\vert\ge\kappa\vert n/d\vert\ge
\kappa\vert n\vert$$
where $d=\gcd(m,n)$; the second
inequality being arrived at as in the proof
above. If $a_1^m=1$, then choose $\kappa=\vert a_{m,e}\vert$,
where $a_{m,e}$ is the coefficient of $z^e$ in
the expansion of $F_m(z)=f^m(z)-z$ and $e\ge2$ is
the smallest integer for which this is non-zero.
Let $b_e$ be the coefficient of the $z^e$ term
of $F_{mn/d}$ ($d=\gcd(m,n)$ as before).
By Lemma~\ref{referee}, $b_e=\frac{n}{d}a_{m,e}$,
so $\vert b_e\vert\ge\kappa\vert n\vert$.
The equation $F_{mn/d}(z)=0$ has
at least $n$ distinct roots of norm $\vert y\vert$
for any point $y$ with least period $n$ under $\phi$;
the rest of the argument proceeds as before.

\begin{pf*}{Proof of Theorem \ref{mainth}}
 We begin the proof with the case where $x$ itself is a periodic point.
 Clearly every periodic point for $\phi$ is also a periodic point
 for a power of $\phi$.
 So it is enough to consider the case where $x$ is a fixed
 point.

 By Proposition \ref{bijective} and Corollary \ref{uni_power}
 the action of $\phi$ in the open unit disk
 $U$ with centre $x$
with respect to the metric $\Delta$ is conjugate to
 the action of a power series $f(z)=\sum_{i=1}^\infty a_iz^i$ on $\cp$.
 If $|a_1|<1$, then it is easy to see that $|f(y)|<|y|$ for any
 $y\in\cp\setminus\{0\}$. So there cannot be any periodic points
 in $\cp$ other than $0$. For $\phi$ this means that there is no
 periodic point $y\in C$ with $\Delta(y,x)<1$.

 Assume now $|a_1|=1$. Then, for some $n$, $|a_1^n-1|<1/p$. As
 before, without loss of generality replace $\phi$ by
 $\phi^n$ and assume $|a_1-1|<1/p$. Lemma
 \ref{referee} shows that there are only finitely many periodic
 points $y\in\cp$ for $f$ with $|y|<r$. This shows the theorem for
 periodic points.

 Assume now $x$ is an arbitrary point. If there is no periodic point $y$
 with $\Delta(x,y)<1$, the statement is trivial. So assume that for some
 periodic point $y$, $R=\Delta(x,y)<1$. By the above
 we know $\Delta(y,y')<r$ holds for only finitely many periodic points $y'$.
 If $R< r<1$, then the ultrametric inequality
 shows that the discs around the centre points $x$ resp.\ $y$ with radius $r$
 agree; the theorem follows.
\end{pf*}

\begin{pf*}{Proof of Theorem~\ref{rational_case}}
The first statement simply specializes Theorem~\ref{mainth}.
For the second, notice that Theorem~\ref{mainth} applied to the
point at infinity implies that $\log^{+}\vert y_n\vert\to0$, so
$\log\vert x-y_n\vert-\log^{+}\vert x\vert\to 0$. The third
follows by factorizing the polynomial $f_n(t)-tg_n(t)$ and noting
that many roots of $f_n(t)-tg_n(t)=0$ are points
with least period $n$ (cf. Remark~\ref{remperpts}).
\end{pf*}

\section{Examples}\label{examplessection}

In this section, we are going to present several examples to
exhibit our main conclusions. The first example explains the earlier
remark that Theorem \ref{rational_case} is a version of Jensen's Formula.

\subsection{Jensen's formula and squaring}\label{jensen}
Assume that $p>2$ and let $f(z)=z^2$. This map gives rise to a
good reduction morphism on $\pri(\Cp)$ with degree 2. Theorem
\ref{rational_case} shows that
$$\lim_{n\rightarrow \infty}2^{-n}\sum_{\zeta^{2^n}=\zeta}
\log \vert \zeta - x\vert _p = \log ^+ \vert x \vert _p.
$$
Working over $\mathbb C$ instead of $\Cp$
the sum on the left would tend to the integral over
the unit circle, and the statement would be exactly
Jensen's Formula. For more details on this point of
view, see \cite{nynts}.

\subsection{Local height on an elliptic curve}\label{duplication}
Let $a$ and $b$ denote elements of $\Cp$
with the
property that $4a^3+27b^2$ does not reduce to zero.
For $p>2$,
the
morphism of degree 4 on $\pri(\Cp)$ defined
by
$$f(X,Y)=\left({X^4-2aX^2Y^2-8bXY^3+a^2Y^4},4{Y(X^3+aXY^2+bY^3)}\right)
$$
has good reduction at $p$.
If the underlying elliptic curve
$$
y^2=x^3+ax+b
$$
is in minimal form at the prime $p$, then our results show that the
(un-normalized)
local height of the point $Q=(x(Q),y(Q))$ can be expressed
as a limit
$$\lim_{n\rightarrow \infty} 4^{-n}\sum_{2^nP=0}\log
\vert x(P)-x(Q)\vert _p.
$$
This example comes about from the duplication map
on the elliptic curve. The reduction condition
guarantees that the elliptic curve has non-singular reduction,
and the condition $p>2$ guarantees that the duplication morphism
has good reduction. This can all be generalized to the multiplication
by $m$ map, and we can also handle the case of an elliptic
curve embedded in projective space in a non-trivial way.

\subsection{Segre embedding}\label{timesmonembeddedcurve}
Let $E$ denote an elliptic curve defined over $\Cp$ with
non-singular reduction. Initially, think of $E$ embedded in
${\mathbb{P}^2}(\Cp)$. For any positive integers $k$ and $\ell$,
map the curve $kE\times\ell E$ to
${\mathbb{P}^8}(\Cp)$ via the Segre embedding (this means
we map $E\to E\times E$ using the map
$P\mapsto(kP,\ell P)\in{\mathbb{P}^2}\times
{\mathbb{P}^2}$ and then embed
the image in $\mathbb{P}^8$ via the Segre
embedding; if $\gcd(k,\ell)=1$ this
is an embedding of $E$). The map $Q\mapsto
mQ$, where $m$ is co-prime to $p$, induces a morphism on this curve
to which Theorem \ref{mainth} applies.

Notice that this is not essentially different
to Section~\ref{duplication}, but gives an
example of how curves can occur in higher-dimensional
projective space.

\section{Background results on the curve and the morphism}\label{Sect_unif}

Let $C\subseteq\pr(\Cp)$ be an irreducible projective smooth curve
with good reduction. For $x\in C$, the ring of regular functions at
$x$ is defined by
\[
  \cO_x=\{f\mid f\mbox{ is a rational regular function on $C$ at }x\}.
\]
For $\pi(x)\in\oC$, define similarly
\[
  \cO_{\pi(x)}=\{f\mid f\mbox{ is a rational
regular function on $C$ at }\pi(x)\}.
\]
Notice that these two rings have quite different properties. For
instance, for $x\in C$, $\cO_x$ is an algebra over $\Cp$, and
\[
 \langle0\rangle\subseteq\{f\in\cO_x\mid f(x)=0\}
\]
is a maximal chain of prime ideals in $\cO_x$. On the other hand,
$\cO_{\pi(x)}$ is an algebra over $O_p$, and
\[
 \langle0\rangle\subseteq\{f\in\cO_{\pi(x)}\mid f(x)=0\}\subseteq\{f\in\cO_{\pi(x)}\mid |f(x)|<1\}
\]
is a maximal chain of prime ideals in $\cO_{\pi(x)}$. Hence, the
Krull dimension of $\cO_x$ is equal to $1$, while that of
$\cO_{\pi(x)}$ is equal to $2$.

\begin{proposition}\label{local}
 Let $x\in C$. Any function
 $z\in\cO_{\pi(x)}$ which vanishes at $x$, and maps modulo $\cp$ to a
 uniformizer
 $\overline{z}$ at $\pi(x)$ for $\oC$, is a uniformizer at $x$ for $C$.
 If $f\in\cO_{\pi(x)}$ vanishes at $x$, then there exists $g\in\cO_{\pi(x)}$
 such that $f=zg$.
 The local power series
 \[
  f(z)=\sum_{i=0}^{\infty}a_iz^i
 \]
 satisfies $a_i\in\Op$.
 Let $y\in C$ with $\Delta(x,y)<1$, then
 \begin{equation}\label{conv}
  f(y)=\sum_{i=0}^{\infty}a_iz(y)^i.
 \end{equation}
\end{proposition}

\begin{pf}
 Let $x=[X]\in C$, $z\in\cO_{\pi(x)}$ be as in the
statement of Proposition \ref{local}.
Assume that $X_0=1$, and work in affine
 coordinates. Let
$$I_0\subseteq\Cp[U_1,\ldots,U_n]\mbox{ and }
J_0=I_0\cap\Op[U_1,\ldots,U_n]$$
be the affine ideals
 corresponding to the homogeneous ideals $I$ and $J$.
Let
 \[
  \cm_x=\langle U_1-X_1,\ldots,U_n-X_n\rangle\subseteq\Cp[U_1,\ldots,U_n]
 \]
be the maximal ideal at $x$, and define a map
 \[
 \theta: \cm_x/ \cm_x^2\rightarrow\Cp^n
 \]
 by
 \[
  \theta(f)=\left(\frac{\partial f}{\partial U_1}(X),\ldots,\frac{\partial
  f}{\partial U_n}(X)\right)\!.
 \]
This is an
 isomorphism between $\cm_x/\cm_x^2$ and $\Cp^n$. Since
 $C$ is smooth with good reduction,
 \begin{equation}\label{n_1}
 \dim (\theta(J_0)+\cp^n)=n-1,
 \end{equation}
where $\cp^n$ as before is the $n$-fold
Cartesian product of $\cp$.
 Let
$$
w_1,\ldots,w_{n-1}\in\theta(J_0)
$$
and $w_n=\theta(z)\in O_p^n$ be elements such that
 $w_1+\cp,\ldots,w_n+\cp$ are linearly independent over
 $\oFp$. Then $w_1,\ldots,w_n\in\Op^n$ are linearly
independent, and the determinant of the matrix formed by
 those vectors is a unit in $\Op$.
Write $v\in\theta(J_0)$ as a linear combination $v=\sum_{i=1}^na_iw_i$ with
 $a_i\in\Op$. Since $w_1,\ldots,w_{n-1},v\in\theta(J_0)$ we
 must have $a_n=0$, for otherwise $\theta(I_0)=\Cp\theta(J_0)$ would be
 $n$-dimensional, contradicting the fact that $C$ is a curve
for which the rank condition~(\ref{rankcondition}) holds.
 This shows that $w_1,\ldots,w_{n-1}$ is a basis for
 $\theta(J_0)$ over $\Op$.
 Therefore $\theta(z)\notin\theta(I_0)$, and $z$ is a uniformizer at
 $x\in C$.

Now define the ideal
 \[
  \cq=\{f\in\cO_{\pi(x)}:f(x)=0\}\subseteq\cO_{\pi(x)}.
 \]
 For any $f\in\cq$, there exists $a\in\Op$ such that
 $\theta(f)-a\theta(z)\in\theta(J_0)$.
 Since $f-az=\frac{h}{g}$ is a rational function which is regular
 at $\pi(x)$ we see that $h(x)=0$ and so $\theta(h)=0$ by the product
 formula for derivatives. It follows that $h\in\cq^2$, so
 \[
  f-az\in\cq^2
 \]
 and therefore
 \[
  \cq=\cO_{\pi(x)}z+\cq^2.
 \]
 Using Nakayama's Lemma~\cite[Th.~2.2]{Matsamura} this shows that
 \begin{equation}\label{Naka}
  \cq=\cO_{\pi(x)}z,
 \end{equation}
 which is the first statement of the proposition.

 For any $f\in\cO_{\pi(x)}$ we can now find $a_0\in\Op$ such that
 $f-a_0\in\cq$. By (\ref{Naka}) there exists $f_1\in\cO_{\pi(x)}$ with
 $f-a_0=f_1z$. For $f_1$ we can find $a_1\in\Op$ and
 $f_2\in\cO_{\pi(x)}$ with $f_1-a_1=f_2z$, and therefore
 \[
 f-(a_0+a_1z)=f_2z^2.
 \]
 Continuing like this gives sequences $a_i\in\Op$ and
 $f_i\in\cO_{\pi(x)}$ such that
 \begin{equation}\label{converge}
  f-\sum_{i=0}^na_iz^i=f_{n+1}z^{n+1}.
 \end{equation}

 Let $y\in C$ with $\Delta(x,y)<1$, then $\pi(x)=\pi(y)$ and
 $|z(y)|=q<1$. Equation (\ref{converge}) shows that
 \[
  \Bigl|f(y)-\sum_{i=0}^na_iz(y)^i\Bigr|\leq q^{-(n+1)},
 \]
 which concludes the proof.
\end{pf}

\begin{proposition}\label{bijective}
 Fix a point $x\in C$.
 There exists a uniformizer $z$ of $C$ at $x$ such that
 $z$ is regular at $\pi(x)$ and $\overline{z}$ is a uniformizer at
 $\pi(x)$.
 The restriction $z:U\to\cp$ of $z$ to $U=\{y\in C\mid
 \Delta(y,x)<1\}$
 is a bijection and its inverse is, in each affine coordinate,
 a convergent power
 series with coefficients in $O_p$. Hence $\Delta(y,y')=|z(y)-z(y')|$ for
 all $y,y'\in U$.
\end{proposition}

\begin{pf}
 Choose $z$ as in Proposition \ref{local}.
 Assume $X_0=1$ and work in the corresponding affine piece. The function $z\in\cO_{\pi(x)}$ which vanishes at $x$
 maps $U$ into $\cp$ and all the affine coordinate projections
 $U_i:C\rightarrow\Cp\cup\{\infty\}$ are elements of
 $\cO_{\pi(x)}$. Applying Proposition \ref{local} concludes the proof.
\end{pf}

Finally, some information about the morphism is needed.

\begin{proposition}\label{phistar}
 Let $C$ be an irreducible projective smooth curve with good reduction.
 Let $\phi:C\rightarrow\prm(\Cp)$ be a morphism of good reduction and
 $x\in C$. Then $\phi$ induces a map
 \begin{eqnarray*}
  &\phi^*:&\cO_{\pi(\phi(x))}\rightarrow\cO_{\pi(x)}\\
      &&  \phi^*(f)=f\circ\phi
 \end{eqnarray*}
\end{proposition}

\begin{pf}
 By the definition of good reduction for
maps, there are forms $F_0,\ldots,F_m$
 such that $\phi$ is represented by $(F_0,\ldots,F_m)$ and for a
 homogeneous coordinate $X$ for $x$
 \[
  |X|=1=|F_0(X),\ldots,F_m(X)|.
 \]
 Let $\frac{F}{G}\in\cO_{\pi(\phi(x))}$ be chosen so that
 the homogeneous coordinate
 $Y=(F_0(X),\ldots,F_m(X))$ of the point $y=\phi(x)$ satisfies $|G(Y)|=1$.
 Then
 \[
  \phi^*z=\frac{F(F_0,\ldots,F_m)}{G(F_0,\ldots,F_m)}=\frac{F^*}{G^*}
 \]
 with $|G^*(X)|=1$, which means that
 $\phi^*\left(\frac{F}{G}\right)\in\cO_{\pi(x)}$.
\end{pf}

Proposition \ref{local} and Proposition \ref{phistar} together
yield the next corollary.

\begin{corollary}\label{uni_power}
 Let $\phi:C\rightarrow C$ be a morphism of good reduction, and let
 $x\in C$ be a fixed point. Then there exists a
 uniformizer $z$, also satisfying Proposition \ref{bijective},
 such that
 \[
  \phi^*z=z\circ\phi=\sum_{i=1}^ {\infty}a_iz^i\mbox{ with
  }a_i\in\Op.
 \]
\end{corollary}

This completes the proof of the tools required for Theorem
\ref{mainth}.

\end{article}

\begin{thebibliography}{1}

\bibitem{beardon-1991}
A. Beardon, \emph{Iteration of {R}ational {F}unctions}, Springer,
New York, 1991.

\bibitem{benedetto}
R. Benedetto, \emph{Reduction, dynamics, and
Julia sets of rational functions},
J. Number Theory, \textbf{86} (2001), 175--195.

\bibitem{call-silverman-1993}
J.~Call and J.H. Silverman, \emph{Canonical heights on
varieties with morphisms},
Compositio Math. \textbf{94} (1993), 163--205.

\bibitem{nynts}
M. Einsiedler, G. Everest and T. Ward, \emph{Morphic heights and periodic
points},
New York Number Theory Seminar, to appear (2002).
{\tt arXiv:math.NT/0204179}

\bibitem{hart}
R. Hartshorne, {\emph Algebraic Geometry},
Springer-Verlag, New York, 1977.

\bibitem{koblitz}
N. Koblitz,  {\emph p-adic Numbers, p-adic {A}nalysis,
and {Z}eta-{F}unctions}, Springer-Verlag,
New York, 1977.

\bibitem{li1}
Hua--Chieh Li, {\emph Counting periodic points of
$p$-adic power series},
Compositio Math. {\textbf{100}} (1996), 351--364.

\bibitem{li2}
Hua--Chieh Li, {\emph p-adic dynamical systems and formal groups},
Compositio Math. {\textbf{104}} (1996), 41--54.

\bibitem{li3}
Hua--Chieh Li, {\emph p-typical dynamical systems and
formal groups},
Compositio Math. {\textbf{130}} (2002), 75--88.

\bibitem{lubin}
J. Lubin, {\emph Nonarchimedean dynamical systems},
Compositio Math. {\textbf{94}} (1993), 321--346.

\bibitem{Matsamura}
H. Matsumura, \emph{Commutative {R}ing {T}heory}, Cambridge
University Press, Cambridge, 1986.

\bibitem{pw}
Y. Puri and T. Ward,
\emph{Arithmetic and growth of periodic orbits},
J. Integer Sequences \textbf{4}:01.2.1 (2001).

\bibitem{thesis}
Juan Rivera Letelier,
\emph{Dynamique des Fonctions Rationelles sur des Corps Locaux},
PhD. thesis, Univ. Paris Sud, Orsay (2000).

\end{thebibliography}
\end{document}